%% file: main.tex
\documentclass{amsart}
\usepackage{local}
\begin{document}

\title{Equivalent computational problems for superspecial abelian surfaces}

\author{Mickaël Montessinos}
\address{ELTE, Faculty of Informatics, H-1053 Budapest, Egyetem tér 1-3, Hungary}
\email{mickael@montessinos.fr}
\thanks{The author is supported by the grant "EXCELLENCE-151343".}



\date{\today}

\begin{abstract}
  We show reductions and equivalences between various problems related to the computation of the endomorphism ring of principally polarised superspecial abelian surfaces. Problems considered are the computation of the Ibukiyama-Katsura-Oort matrix and computation of unpolarised isomoprhisms between superspecial abelian surfaces.
\end{abstract}

\maketitle


\input{sections/introduction}
\input{sections/preliminaries}
\input{sections/problems}
\input{sections/endomorphism_rings}
\input{sections/isomorphism}
\input{sections/endo_to_iko}

\bibliographystyle{plain}
\bibliography{biblio}
\end{document}

%% file: sections/introduction.tex
\section{Introduction} The problem of computing the endomorphism ring of a supersingular elliptic over a finite field has been well studied in recent years. It was proven equivalent to the Isogeny Path Problem and the One Isogeny Problem \cite{herledan2025unconditional}. In particular, the hardness of the endomorphism ring computation underlies the security of most one-dimensional isogeny-based cryptographic protocols.

The precise statement of the problem depends on the choice of representation for the endomorphism ring of a supersingular elliptic curve. Essentially, one may either wish to know an efficient algorithm to evaluate \(4\) endomorphisms which form a \(\ZZ\)-basis of the endomorphism ring, or to know the basis of a maximal order in the quaternion \(\QQ\)-algebra \(B_{p, \infty}\) which ramifies at the prime \(p\), the characteristic of the base field of the curve. The two forms of the problem were proved to be equivalent, first heuristically \cite{eisentrager2018supersingular}, then unconditionally under GRH \cite{wesolowski2022supersingular} and finally unconditionally \cite{herledan2025unconditional}.

The main tool appearing in reductions from most problems underlying the security of isogeny-based schemes is the so-called \(\mathsf{KLPT}\) algorithm. Its original version was published in \cite{kohel2014quaternion}, and relied on specific conditions regarding the characteristic of the base field of the elliptic curves involved, as well as some heuristic considerations. Several subsequent works improved both the performance and the assumptions of the algorithm, and a general version relying only on GRH was given in \cite{wesolowski2022supersingular}. This algorithm uses quaternion arithmetic to compute an isogeny of power-smooth degree between elliptic curves with known endomorphism rings. We note however that the more recent results from \cite{herledan2025unconditional} overcome the reliance on GRH by circumnavigating the need for using KLPT.

In recent years, applications of the algorithmic theory of principally polarised superspecial abelian varieties of dimension greater than one has emerged, both regarding the construction of cryptographic schemes \cite{castryck2020hash} and the cryptanalysis of schemes using supersingular elliptic curves \cite{castryck2023efficient,maino2023direct,robert2023breaking}. As in the case of one dimensional protocols, most security assumptions rely on the difficulty of computing polarised isogenies between principally polarised superspecial abelian surfaces (PPSAS). As for maximal quaternion orders, a \(\mathsf{KLPT}^2\) algorithm working over hermitian matrices in \(M_2(B_{p, \infty})\) and allowing for the computation of polarised isogenies of power-smooth degrees between surfaces was given in \cite{castryck2025klpt}.

In the higher dimensional case, since all unpolarised superspecial abelian varieties are isomorphic, their endomorphisms rings are as well. In fact, they are all isomorphic to the ring \(M_2(O)\), where \(O\) is any maximal order in the quaternion algebra \(B_{p, \infty}\). The input to the \(\mathsf{KLPT}^2\) algorithm is therefore not a characterisation of the endomorphism ring, but rather the data of a certain Hermitian matrix in \(M_2(O)\), which we call the Ibukiyama-Katsura-Oort matrix, for reasons that will be made clear in \Cref{sec:prelim}.

Much like in the one dimensional case, a question naturally arises of what the various alternative computational representations of a usable "known endomorphism ring" are, and of their equivalences. The knowledge the Ibukiyama-Katsura-Oort matrix of a principally polarised superspecial abelian surface is analogous to the knowledge of a maximal quaternion order isomorphic to the endomorphism ring of a curve, in the sense that it is the input to the relevant \(\mathsf{KLPT}^{2}\) algorithm. Another possible representation is the knowledge of efficient algorithms representing \(16\) endomorphisms of an abelian variety \(A\) which form a basis of the endomorphism ring. If the variety \(A\) is represented as either a product of elliptic curves or the Jacobian of a hyperelliptic curve of genus \(2\), one may hope that an efficient representation of these endomorphisms be enough to compute their image by the Rosati involution related to the principal polarisation implied by the representation of \(A\).

Another problem to consider is that of the computation of an unpolarised isomoprhism between two PPSAS. It is known that the restriction of this problem to PPSAS presented as products of supersingular elliptic curves may be solved when the endomorphism rings of the curves are known \cite{gaudry2025computing}. This problem is also a natural candidate for equivalence to the computation of the Ibukiyama-Katsura-Oort matrix of a PPSAS: one may recover the endomorphisms of the factors of a product surface from the Ibukiyama-Katsura-Oort matrix, and then compute isomorphisms to other products. Conversely, computing an unpolarised isomorphism from a variety \(A\) to a reference surface \(A_0\) may allow to recover the Ibukiyama-Katsura-Oort matrix.

Finally, one may fix a supersingular curve \(E_0\) with known endomorphism ring \(O_0\), and consider the structure of the Hermitian left \(O_0\)-module \(\Hom(A, E_0\) associated to a PPSAS \cite{robert2024module}. We do not consider this problem in this work, but it is the object of \cite{page2026conversion}.

\subsection{Our contributions}
Let \(A\) be a PPSAS. If \(A\) is presented as a product of supersingular elliptic curves, we prove that the following problems are equivalent under GRH and the assumptions necessary to use the \(\mathsf{KLPT}^2\) algorithm:
  \begin{itemize}
    \item Computing the Ibukiyama-Katsura-Oort matrix of \(A\);
    \item Computing an effective representation of endomorphisms of \(A\) that form a basis of \(\End(A)\) as a \(\ZZ\)-module;
    \item Computing an unpolarised isomorphism from \(E_0^2\) (or any product of curves with known endomorphism rings) to \(A\)
    \end{itemize}

  If \(A\) is presented as the Jacobian surface of a genus \(2\) hyperelliptic curve instead, we prove that the problem of computing the Ibukiyama-Kastura-Oort matrix of \(A\) is equivalent to the problem of computing a so-called good representation of its endomorphism ring, which contains effective representations of endomorphisms of \(A\) as above, but also some additional information. We also prove that the problems of computing an effective representation of the endomorphism ring of \(A\) and computing an unpolarised isomorphism to any other PPSAS with known Ibukiyama-Katsura-Oort matrix both reduce to the computation of the Ibukiyama-Katsura-Oort matrix of \(A\).

  Theses results are stated formally as \Cref{thm:main}, and a diagram presenting the various reductions is given there.

  It should be noted that, albeit polynomial, the reduction from computing the Ibukiyama-Katsura-Oort matrix of a variety from its endomorphism ring with additional information is not practical. Indeed, the algorithm for computing an isomorphism between the endomorphism algebra \(\End(A)
  otimes \QQ\) and \(M_2(B_{p, \infty})\) has a very large constant part in its complexity.

  We provide an implementation of the algorithm from \Cref{lemma:iso_comput_2_2} at the following address:

  \url{https://gitlab.com/Mickanos/isomorphisms-to-superspecial-jacobians}

\subsection{Methods}

The two main technical contributions of this work are the reduction from computing unpolarised isomorphisms to computing Ibukiyama-Katsura-Oort matrices, and from computing Ibukiyama-Katsura-Oort matrices to computing the endomorphism ring of a PPSAS and some additional information.

The problem of computing unpolarised isomorphisms between products of supersingular elliptic curves knowing their endomorphism rings being solved, it is enough to provide an algorithm which, given a superspecial Jacobian surface and its Ibukiyama-Katsura-Oort matrix, compute an isomorphism to any product of elliptic curves, together with their endomorphism rings.

This is done by first computing a chain of \((2,2)\)-isogenies from a square \(E_0^2\), where \(E_0\) is a curve with known endomorphism ring. Then, given the kernel \(K\) of a \((2,2)\) isogeny \(E_1 \times E_2 \to J(C)\), one computes an unpolarised automorphism of \(E_1 \times E_2\) which sends \(K\) to a group \(K'\) which is the kernel of a pair of parallel \(2\)-isogenies going respectively from \(E_1\) and \(E_2\) to some curves \(E_3\) and \(E_4\). One then recovers an isomorphism from \(E_3 \times E_4\) to \(J(C)\) by completing a square diagram.

Considering a PPSAS \(A\) with a known basis of endomorphisms, and the data of the action of these endomorphisms on the space of global differentials, as well as the matrix of the Rosati involution with respect to that basis, a natural strategy for recovering the Ibukiyama-Katsura-Oort matrix is to compute an isomorphism \(\theta: \End(A) \otimes \QQ \to \End(E_0^2) \otimes \QQ\), where \(E_0\) is the base curve of the Ibukiyama-Katsura-Oort correspondence, and then carry the Rosati involution through that isomorphism. One then finds an element \(\gamma \in \End(E_0^2)\) such that \(\gamma \iota(\End(A)) \gamma^{-2} = \End(E_0^2)\) and carries over the Rosati involution of \(\End(A)\) first through \(\iota\) and then conjugates it by \(\gamma\) as well.

This strategy fails to account for the existence of automorphisms of \(\End(A)\) which are not induced by an automorphism of \(A\), as then the isomorphism that was constructed from \(\End(A)\) to \(\End(E_0^2)\) may itself not be induced by an unpolarise disomorphism from \(A\) to \(E_0^2\). Whether a given isomorphism of endomorphism rings is as above can be detected by checking the effect of endomorphisms on global differentials, and then corrected by post-composing another such "faulty" automorphism of \(E_0^2\) if necessary.

\subsection{Concurrent works}
\begin{itemize}
  \item The problem of computing unpolarised isomorphisms between products of supersingular elliptic curves is treated in \cite{gaudry2025computing}. It is proved that this problem reduces to computing the endomorphism rings of all curves involved.

  \item Computations equivalent to the computation of the Ibukiyama-Katsura-Oort matrix of a PPSAS are discussed in \cite{page2026conversion}, where more emphasis is put on the Hermitian module \(\Hom(A, E_0)\) associated to a PPSAS \(A\). In contrast, this last work does not show an equivalence between the computation of the Ibukiyama-Katsura-Oort matrix of a PPSAS and the computation of its endomorphism ring, eithe in good or effective representation.
\end{itemize}

\subsection{Structure of this work}

The remainder of this article is organised in the following way:
\begin{itemize}
  \item \Cref{sec:prelim} recalls the definitions of the main concepts that we use, states some known results and show a few general lemmas necessary to our results.
  \item \Cref{sec:problems} states the various computational problems that we consider, as well as our main result, \Cref{thm:main}, which gives the relationship between these problems.
  \item \Cref{sec:endo} shows reduction of the computation of the endomorphism ring of a PPSAS to various problems.
  \item \Cref{sec:iso} completes the work of \cite{gaudry2025computing} in showing that the computation of unpolarised isomorphisms between PPSAS reduces to the computation of their Ibukiyama-Katsura-Oort matrices.
  \item \Cref{sec:endo_to_iko} shows that the computation of the Ibukiyama-Katsura-Oort matrix reduces to the computation of a good representation of the ring of endomorphisms.
\end{itemize}

\subsection{Open problems}

The situation we paint in this work is somewhat neater when we restrict to products of supersingular elliptic curves. This is due to the existence of strong results for supersingular elliptic curves, which we are able to take advantage of when working with a product. In contrast, we do not prove that all the computational problems considered are equivalent when Jacobian surfaces are included.

We leave to future work the consideration of how more appropriate representations of maps from a Jacobian variety may allow to complete our equivalence diagram.

\subsection{Notations}
We let \(p\) be an odd prime number and \(q = p^2\), and we fix the finite field \(k = \FF_q\). Unless stated otherwise, every variety that we consider is defined over \(k\). We assume that \(p\) is so that the \(\mathsf{KLPT}^2\) algorithm may be used. The current state of the art exposed in \cite{castryck2025klpt} is stated with the hypothesis that \(p = 3 \mod 4\), although it is suggested that any \(p\) such that there is a maximal order \(O_0\) of the quaternion algebra \(B_{p,\infty}\) with discriminant \(p\) which contains a quadratic order with very small discriminant.

We also adopt some notations related to the context of \(\mathsf{KLPT}^2\): we fix a supersingular elliptic curve \(E_0/k\) with known endomorphism ring, for instance \(E_0: y^2 = x^3 + x\) if it is assumed that \(p = 3 \mod 4\). Then, we let \(O_0\) be the usual maximal order of \(B_{p, \infty}\) isomorphic to \(\End(O_0)\). We set \(A_0 = E_0^2\), and identify the maximal order \(M_2(O_0)\) of \(M_2(B_{p,\infty})\) in the natural way from the identification of \(O_0\) as \(\End(E_0)\).

If \(M \in M_2(B_{p, \infty})\), we write \(M^*\) for the conjugate transpose, where the conjugation is the special involution of the quaternion algebra \(B_{p, \infty}\).

Throughout this work, we refer to "the assumptions for applying the \(\mathsf{KLPT}^2\) algorithm". This is phrased in anticipation of future works providing a relaxation of the assumptions and heuristics from \cite{castryck2025klpt}. We state the assumption of the generalised Riemann hypothesis separately, as it is also a condition for applying the \(\mathsf{KLPT}\) algorithm, although it seems unlikely that an unconditional version of \(\mathsf{KLPT}^2\) would appear before an undoncitional version of \(\mathsf{KLPT}\).

\paragraph{Acknowledgments}
The author is grateful to Péter Kutas for pointing out the manuscript \cite{page2026conversion} and to the authors of this manuscript, Aurel Page, Damien Robert and Julien Soumier for sharing their manuscript. The author is also grateful to Marc Houben for an insightful conversation about \Cref{lemma:preimage_2_2_isogeny}.

%% file: sections/preliminaries.tex
\section{Preliminaries}\label{sec:prelim}

\subsection{Principally polarised superspecial abelian surfaces}
    \subsubsection{Polarisations and the Rosati involution}

    We briefly recall here the basic definitions and results of the theory of polarisations of abelian varieties. We direct the reader to \cite{edixhovenabelian,milne1984abelian} for a more detailed exposition.

    Let \(A\) be an abelian variety, then there exists a variety \(\widehat{A}\) \(\widehat{A}(k) \simeq \Pic^0(A)\). An isogeny \(f \colon A \to A'\) admits a dual isogeny \(\widehat{A'} \to \widehat{A}\) induced by the pullback operation \(f^*\) on divisors of \(A'\). 

    \begin{definition}
        A \emph{polarisation} on \(A\) is an isogeny \(\lambda\colon A \to \widehat{A}\) for which there exists an ample divisor \(D\) on \(A_{\overline{k}}\) such that \(\lambda(P) = [t_{-P}(D) - D]\).

        The \emph{degree} of a polarisation is its degree as an isogeny. A \emph{principal polarisation} is a polarisation of degree \(1\), and therefore an isomorphism.

        A \emph{(principally) polarised abelian variety} is the data \((A, \lambda)\) of an abelian variety \(A\) together with a (principal) polarisation \(\lambda\) on \(A\).
    \end{definition}

    We next define the condition for an isogeny to be a morphism between polarised abelian varieties:

    \begin{definition}
        Let \((A, \lambda)\) and \((A', \lambda')\) be polarised abelian varieties. An isogeny \(f\colon A \to A'\) is \emph{polarised} if the following equation holds for some \(N \in \ZZ\):
        \[\widehat{\varphi} \circ \lambda' \circ \varphi = [N]_{\widehat{A}} \circ \lambda.\]

        Two polarisations \(\lambda, \lambda'\) on an abelian variety \(A\) are \emph{equivalent} if there is a polarised isomorphism from \((A, \lambda)\) to \((A, \lambda')\).
    \end{definition}

    The choice of a polarisation on \(A\) is reflected on the ring \(\End(A)\) by the corresponding Rosati involution. While the Rosati involution is defined over \(\End(A) \otimes \QQ\) for any polarisation, we shall only need the notion for principal polarisations:

    \begin{definition}
        Let \((A, \lambda_A)\) and \((B, \lambda_B)\) be principally polarised abelian varieties. Let \(\varphi\colon A \to B\) be a homomorphism of algebraic groups. Then, we set
        \[\begin{array}{rccc}
          r_{\lambda_A,\lambda_B}\colon &\Hom(A, B) &\to &\Hom(A, B) \\
            & f &\mapsto & \lambda_A^{-1} \circ \widehat{f} \circ \lambda_B.
        \end{array}\] 
        In particular, if \(A = B\) and \(\lambda_A = \lambda_B = \lambda\), we write \(r_\lambda = r_{\lambda_A,\lambda_B}\). In this case, the map \(r_\lambda\) is the well known \emph{Rosati involution} on \(\End(A)\).
    \end{definition}

    If \((A, \lambda)\) is a principally polarised abelian variety and the polarisation is clear from contex, we also write \(r_A\) for the Rosati involution on \(\End(A)\).
    Then, it is well known that and a fundamental fact that the bilinear map
    \[\begin{array}{ccc}
      \End(A) \times \End(A) &\to &\ZZ\\
      (f, g) &\mapsto &\langle f, g\rangle \coloneqq \Tr(f \circ r_A(g))
    \end{array}\]
    is positive definite. \cite[Theorem 17.3]{milne1984abelian}

    \begin{example}\label{ex:natural_principal_representations}
    We give some examples of abelian varieties with canonical principal polarisations:
        \begin{itemize}
            \item An elliptic curve \(E\) has a canonical principal polarisation with ample divisor \([O]\), the (arbitrarily chosen) neutral element. Concretely, it is given by the identification \(P \mapsto [P - O]\), which is often used to define the group law on \(E\). The corresponding Rosati involution then simply sends an endomorphism \(f\) to its traditionally defined dual \(\widehat{f}\).
            \item A product of elliptic curves \(E_1 \times \dots \times E_n\) admits a canonical polarisation, which is the product of the canonical polarisations of the \(E_i\). Its ample divisor is then \(\sum_{i = 1}^n E_1 \times \dots \times E_{i-1} \times \{\infty_i\} \times E_{i + 1} \times \dots \times E_n\).
            \item The Jacobian variety of an algebraic curve admits a natural principal polarisation. If \(C\) is an algebraic curve of genus \(2\) and \(\infty \in C\) is a fixed point, there is an injective map \(u\colon C \to J(C)\) sending a point \(P\) to \([P - \infty]\). Then, \(u(C)\) is an ample divisor of \(J(C)\) which defines the canonical polarisation.
        \end{itemize}
    \end{example}

    \subsubsection{Superspecial abelian varieties}

    Superspecial abelian varieties come as one of the possible generalisations of supersingular elliptic curves to larger dimension.

    \begin{definition}
        An abelian variety \(A\) is superspecial if there is a supersingular elliptic curve \(E\) and \(g \in \NN\) such that \(A_{\overline{k}} \simeq E^g_{\overline{k}}\).
    \end{definition}

    A well known result by Deligne, Ogus and Shioda states that all superspecial abelian variety of equal dimension are isomorphic over \(k\).

    In fact, we may choose our models so that an even stronger statement is true:
    \begin{definition}\label{def:maximal_variety}
      An abelian variety is called \emph{maximal} if \(\#A(k) = (p + 1)^{2\dim A}\).
    \end{definition}

    \begin{proposition}\cite[Section 5]{jordan2018abelian}\label{lemma:maximal_varieties}
      Let \(A\) be a maximal abelian variety. Then, we have the following:
      \begin{enumerate}
        \item \(A\) is superspecial;
        \item If \(E\) is a maximal elliptic curve, there is a \(k\)-rational isomorphism from \(A\) to \(E^{\dim A}\);
        \item If \(B\) is another maximal abelian variety of same dimension, \(A\) is isomorphic to \(B\), and any homomorphism \(A_{\overline{k}} \to B_{\overline{k}}\) is the base extension of a homomorphism \(A \to B\).
      \end{enumerate}
    Furthermore, every superspecial abelian variety defined over \(\overline{k}\) admits a maximal model defined over \(k\).
    \end{proposition}

    We state a well known classifying result on principally polarised superspecial abelian surfaces.
    \begin{lemma}
      Let \(A\) be a PPSAS. Then there is a polarised isomorphism from \(A\) to either the Jacobian of a genus \(2\) hyperelliptic curve or the product of two supersingular elliptic curve.
      Furthermore, two Jacobians \(J(C)\) and \(J(C')\) are isomorphic if and only if the curves \(C\) and \(C'\) are isomorphic, and two products are isomorphic if and only if their factors are so up to reordering.
    \end{lemma}

    \begin{proof}
      The first result follows from \cite[Theorem 1]{flynn2019genus} and \Cref{lemma:maximal_varieties}. It follows from the criterion of Matsusaka-Ran \cite[Theorem 14.46]{edixhovenabelian} that a PPSAS may not be both isomorphic to a Jacobian of a genus \(2\) hyperelliptic curve and to a product of supersingular elliptic curves, and that two such products are isomorphic if and only if their factors are isomorphic up to reordering. The statement on isomorphisms of Jacobians is the Torelli theorem \cite[Corollary 12.2]{milne1984jacobian}. 
    \end{proof}

    \subsubsection{The Ibukiyama-Katsura-Oort classification}
    Ibukiyama, Katsura and Oort introduced a bijection between the set of polarisations on the abelian variety \(E_0^2\) and a certain set of matrices lying in \(M_2(O_0)\) \cite{ibukiyama1986supersingular}. We recall their main definitions and results, although we mainly follow the notations and exposition from \cite[Section 2.3]{castryck2025klpt}.

    Since all superspecial abelian surfaces are isomorphic to one another, classifying principal polarisations over \(A_0 = E_0^2\) up to equivalence is the same as classifying principally polarised superspecial abelian surfaces up to polarised isomorphism. 

    We let \(\lambda_0\) be the natural product polarisation on \(A_0\). Then, if \(\lambda\) is a principal polarisation on \(A_0\), we set
    \[\mu(\lambda) = \lambda_0^{-1} \lambda \in \End(A_0) = M_2(O_0).\]
    We then have the following result.
    \begin{proposition}[{\cite[Theorem 2.7]{castryck2025klpt}}]
        The map \(\mu\) is injective, and its image in \(M_2(O_0)\) (using the identification \(\End(E_0) \simeq O_0\)) is
        \[\Mat(A_0) \coloneqq \left\{\begin{pmatrix} s & r \\ \overline{r} & t\end{pmatrix},\quad s, t \in \NN, r \in O_0, st - r\overline{r} = 1 \right\}\]
    \end{proposition}

    In several formulas, the principal polarisation \(\lambda\) may be replaced almost verbatim with \(g = \mu(\lambda)\) in order to reframe computations within \(M_2(O_0)\).

    We cite two such cases: 
    \begin{lemma}\label{lemma:ros_iso_IKO}
        Let \(\lambda\) be a principal polarisation on \(O_0\), and let \(\gamma \in M_2(O_0)\) represent an endomorphism of \(A_0\). We then have
        \[r_\lambda(\gamma) = \mu(\lambda)^{-1} \gamma^* \mu(\lambda).\]

        Let \(\lambda, \lambda'\) be principal polarisations, and let \(\gamma \in M_2(O_0)\) be an endomorphism. Then \(\gamma\) is a polarised isogeny from \((A_0, \lambda)\) to \((A_0, \gamma')\) if and only if the following equation holds for some \(N \in \NN\):
        \[\gamma^* \mu(\lambda') \gamma = N \mu(\lambda).\]
    \end{lemma}

    In particular, it is easy to recover the matrix \(\mu(\lambda)\) from the data of the Rosati involution \(r_\lambda\):
    \begin{lemma}\label{lemma:compute_mu_from_rosati}
        There is a polynomial-time algorithm which, given a matrix for the Rosati involution \(r_\lambda\) with respect to some principal polarisation \(\lambda\) of \(A_0\), computes \(\mu(\lambda)\).
    \end{lemma}

    \begin{proof}
        Let \(B\) be a basis of \(M_2(O_0)\). Then, \(\mu(\lambda)\) is the solution in \(\Mat(A_0)\) of the following system of linear equations (in variable \(g\)):
        \[g r_\lambda(\gamma) - \gamma^* g = 0,\quad \forall \gamma \in B.\]
        Such a system can be solved in polynomial-time in \(\QQ \Mat(A_0)\), and then a solution may be rescaled to lie in \(\Mat(A_0)\).
    \end{proof}

    \begin{remark}
        It follows directly from \Cref{lemma:ros_iso_IKO} that two principal polarisations \(\lambda\) and \(\lambda'\) on \(A_0\) are equivalent if and only if there exists \(\gamma \in GL_2(O_0)\) such that
            \[\gamma^* \mu(\lambda') \gamma = \mu(\lambda).\]
        This defines an equivalence relationship on \(\Mat(A_0)\). In particular, if \((A, \lambda_A)\) is any PPSAS, it is isomorphic to \((A_0, \lambda)\) for some principal polarisation \(\lambda\), and we write \(\mu(A_0)\) for \(\mu(\lambda)\). It is then understood that \(\mu(A_0)\) is only defined up to equivalence in \(\Mat(A_0)\).
        We call \(\mu(A_0)\) the \emph{Ibukiyama-Katsura-Oort matrix of \(A\)}.
    \end{remark}

\subsection{Computational Representations}
    \subsubsection{Principally polarised superspecial abelian surfaces}
        While it is possible to give equations for superspecial abelian surfaces, it is generally preferable to represent them as either products of supersingular elliptic curves or Jacobians of hyperelliptic curves of genus \(2\). In both cases, the isomorphy class of the principally polarised abelian surface is precisely described by the isomorphy classes of the curves that induce it.

        Recall that the isomorphy class of an elliptic curve \(E\) is described by its \(j\)-invariant, \(j(E) \in k\). Likewise, the isomorphy class of a hyperelliptic curve of genus 2 is described by its Igusa invariant \cite{igusa1960arithmetic}.

        In both cases, knowing the invariants allows us to provide a model for the curve, and compute arithmetical operations. In general, we assume that the models used in this work are all maximal, so that every homomorphism is defined over \(k\).

        In the sequel, any PPSAS is computationally represented by either an unordered pair of \(j\)-invariants or by a projective vector of Igusa invariants.

    \subsubsection{Isogenies between abelian surfaces}
        Since principally polarised abelian surfaces are not directly represented as embedded projective varieties, isogenies may not simply be represented using rational functions.

        We follow the definitions from \cite{wesolowski2022supersingular,anni2025computation} and introduce the following definitions:
        \begin{definition}
        Let \(\varphi\colon A \to B\) be a homomorphism between two abelian varieties defined over a finite field \(\FF_q\). An \emph{efficient representation} of \(\varphi\) with respect to a given algorithm is some data \(D_\varphi \in \{0,1\}^*\) such that, on input \(D_\varphi\) and \(P \in A(\FF_q)\), the algorithm returns the evaluation \(\varphi(P)\) in polynomial time in \(\mathrm{length}(D_\varphi)\) and \(\log q\).
      \end{definition}

      When \(A\) and \(B\) are principally polarised abelian varieties, some additional information may be required in order to relate the isogeny to the polarisations of \(A\) and \(B\). For this reason, we also use the definition of a good representation of an isogeny, as is introduced in \cite{anni2025computation}. We add one extra piece of information to be carried, about the behaviour of the isogeny over differentials of \(B\). The reason for this addition will be made clear in \Cref{sec:endo_to_iko}.

      \begin{definition}
        Let \(A, \lambda_A\) and \(\lambda_B\) be principally polarised abelian varieties defined over a finite fild \(\FF_q\). Let \(\varphi\colon A \to B\) be a homomorphism. A \emph{good representation} of \(\varphi\) is a tuple \(\theta, \theta^\dagger, D, \theta^*\), where \(\theta\) is an efficient representation of \(\varphi\), \(\theta^\dagger\) is an efficient representation of \(r_{\lambda_A, \lambda_B}(\varphi)\), \(D\) is an integer such that \(\Tr(\varphi \circ r_{\lambda_A, \lambda_B}(\varphi)) \leq D\) and \(\theta^*\) is a matrix representing the \(\FF_q\)-linear map from \(H^0(B, \Omega_B)\) to \(H^0(A, \Omega_A)\) by \(\varphi\) with respect to any bases of these vector spaces.
      \end{definition}

        We provide two examples, which will be omnipresent in the sequel:

        \begin{example}\label{ex:isogeny_matrix}
          Let \(E_1, E_2, E'_1, E'_2\) be supersingular elliptic curves and let \(\varphi\colon E_1 \times E_2 \to E'_1 \times E'_2\) be a homomorphism. This homomorphism decomposes as a matrix \((\varphi_{ij})_{1 \leq i,j \leq 2}\), where \(\varphi_{ij} \in \Hom(E_j, E'_i)\). Then, the data of efficient representations \(\theta_{ij}\) respectively of the \(\varphi_{ij}\) is an efficient representation of \(\varphi\).

          Furthermore, letting \(\lambda\) and \(\lambda'\) be the product polarisations respectively on \(E_1 \times E_2\) and \(E'_1 \times E'_2\), then as discussed in \cite[Section 2.2]{castryck2025klpt}, we have
        \[r_{\lambda,\lambda'}(\varphi) = (\widehat{\varphi}_{ji}).\]

        Since differentials are compatible with products of varieties, we also have that \(\varphi^* = (\varphi^*_{ij})\) as a block matrix.
        \end{example}

        \begin{example}[(2,2)-isogenies]
          A \((2,2)\)-isogeny is an isogeny between abelian surfaces, whose kernel is maximally isotropic with respect to the Weil pairing, and therefore isomorphic to \(\ZZ/2\ZZ \times \ZZ/2\ZZ\). If \(A\) is a principally polarised abelian surfaces, there exists formula which, given such a subgroup \(K\), yield a polarised \((2,2)\)-isogeny as well as its codomain. There are also explicit criteria to decide whether the codomain is a Jacobian variety or a product of elliptic curves. The explicit formulas and criteria are given in \cite[Section 3.2]{castryck2020hash} and \cite[Proposition 4 and its proof]{howe2000large}.
            If \((A, \lambda_A)\) and \((B, \lambda_B)\) are principally polarised abelian surfaces, and \(K \subset A\) is a subgroup isomorphic to \(\ZZ/2\ZZ \times \ZZ/2\ZZ\), and \(B\) is the codomain of the \((2,2)\) isogny with kernel \(K\), then \(K\) is a good representation of the isogeny with kernel \(K\).
            
            For our purposes, we note that if \(A = E_1 \times E_2\) is a product of elliptic curves, and the subgroup \(K \subset A[2]\) admits generators of the form \((P, 0)\) and \((0, Q)\), then the codomain of the polarised \((2,2)\)-isogeny with kernel \(K\) is the product \(E_1 / P \times E_2 / Q\).
        \end{example}

        \begin{remark}
          Not all subgroups \(K \subset A\) that are isomorphic to \(\ZZ/2\ZZ \times \ZZ/2\ZZ\) are the kernel of an isogeny. Indeed, if \(A = E_1 \times E_2\) is a product of elliptic curves, the subgroup \(E_1[2] \times \{\infty_2\}\), for instance, is not maximally isotropic with respect to the Weil pairing.
        \end{remark}

    \subsubsection{Endomorphism rings}
        Throughout this work, we regularly use the phrase "known endomorphism ring" about elliptic curves and principally polarised abelian surfaces.

        In order to make this statement precise, we introduce the following definition:
        \begin{definition}
          Let \(A\) be a abelian variety. 
          By an \emph{efficient representation} (resp. a \emph{good representation}) of the endomorphism ring of \(A\), we mean efficient representations (resp. good representations) of endoomrphisms of \(A\) that form a basis of \(\End(A)\) as a \(\ZZ\)-module.
          We also say that the endomorphism ring of \(A\) is \emph{weakly known} (resp. \emph{known}) if we are given efficient representations (resp. good representations) of endomoprhism of \(A\) forming a basis of \(\End(A)\) as a \(\ZZ\)-module.
        \end{definition}

        It is clear that if the endomorphism ring of \(A\) is known, then it is weakly known. In the case of a supersingular elliptic curve \(E\), it is known that, under GRH, there is an equivalence between weakly knowing and knowing the endomorphism ring.

        \begin{lemma}[{\cite[Theorems 8.1 and 8.3]{wesolowski2022supersingular}}]\label{lemma:curve_knowing_endo}
            Under GRH, there is a polynomial-time algorithm which takes as input a supersingular elliptic curve \(E\) with weakly known endomorphism ring and outputs a good representation of \(\End(E)\).
        \end{lemma}

        We introduce a notation for the isomorphism of endomorphism rings induced by an unpolarised isomorphism of abelian varieties:
        \begin{definition}
          Let \(\varphi\colon A_1 \to A_2\) be an unpolarised isomorphism of abelian varities. We then define the following isomorphism of rings:
          \[\begin{array}{rccc}
            \varphi_*\colon & \End(A_1) &\to &\End(A_2) \\
            &f &\mapsto & \varphi \circ f \circ \varphi^{-1}
          \end{array}\]
        \end{definition}

\subsection{Some subroutines}

\begin{lemma}\label{lemma:invert_isom_product}
    Let \(E_1, E_2, E_3, E_4\) be superspecial elliptic curves, and let \(\varphi\) be an isomorphism from \(E_1 \times E_2 \to E_3 \times E_4\) represented as in \Cref{ex:isogeny_matrix}.
    Then the isomorphism \(\varphi^{-1}\) may be computed in polynomial time.
\end{lemma}

\begin{proof}
    By \cite[Corollary 4.11]{gaudry2025computing} and its proof, we may compute an isomorphism \(\tilde{\varphi}\colon E_3 \times E_4 \to E_1 \times E_2\) such that \(\varphi \circ \tilde{\varphi}\) is an automorphism of the form
    \[\varphi \circ \tilde{\varphi} = \begin{pmatrix} a & \widehat{\psi} \\ \psi & d \end{pmatrix},\]
    where \(a, d \in \ZZ\) and \(\psi \colon E_3 \to E_4\) is an isogeny.
    Then, by \cite[Remark 5.1]{gaudry2025computing}, we have
    \[(\varphi \circ \tilde{\varphi})^{-1} = \begin{pmatrix} d & -\widehat{\psi} \\ -\psi & a \end{pmatrix}.\]
    It follows that \(\varphi^{-1} = \tilde{\varphi} \circ (\varphi \circ \tilde{\varphi})^{-1}\) may be computed in polynomial time.
\end{proof}

\begin{lemma}\label{lemma:preimage_2_2_isogeny}
    Let \(E_1, E_2\) be elliptic curves and let \(A\) be an abelian surface. Let \(f\colon E_1 \times E_2 \to A\) be a \((2,2)\)-isogeny, and let \(P \in A\). Then a point \(P \in E_1 \times E_2\) such that \(f(P) = Q\) may be computed in polynomial time.
\end{lemma}

\begin{proof}
  Let \(K = \Ker(f)\), and let \(K' = f((E_1 \times E_2)[2])\). Let \(g\) be the \((2,2)\)-isogeny with kernel \(L\) and codomain \(E_1 \times E_2\). Then we have \(f \circ g = [2]_A\).
  For a point \(Q \in A\), compute \(Q' \in A\) such that \(2Q' = Q\), and then set \(P = g(Q')\).
\end{proof}

\begin{lemma}\label{lemma:klpt2_2_2}
    Let \(A_1, A_2\) be principally polarised superspecial abelian varieties, with known Ibukiyama-Katsura-Oort matrices \(\mu(A_1)\) and \(\mu(A_2)\). Then, assuming the hypotheses for the KLPT2 algorithm, we may compute in polynomial time a chain \((f_i)_{0 \leq i \leq e}\) of \((2,2)\) isogenies \(f_i \colon A'_i \to A'_{i+1}\) such that \(A'_0 = A_1\) and \(A'_{e+1} = A_2\).
\end{lemma}

\begin{proof}
    This is a combination of Algorithms \(\mathsf{KLPT}^2\) and \(\mathsf{MatrixToIsogeny22}\) from \cite{castryck2025klpt}.
\end{proof}

\begin{lemma}\label{lemma:explicit_isom_CSA}
    There exists a polynomial-time algorithm which, given structure constants for two central simple \(\QQ\)-algebras of degree bounded by \(4\) and a basis of a maximal order in each of them, computes an isomorphism between these algebras.
\end{lemma}

\begin{proof}
    This is stated for quaternion algebras in \cite[Proposition 4.1]{csahok2022explicit}, but the argument adapts readily to central simple algebras of larger, albeit bounded, degree.
\end{proof}

\begin{lemma}\label{lemma:structure_constants}
  If the endomorphism ring of \(A\) is known with basis \(B\), the structure constants of \(\End(A)\) with respect to its basis \(B\) may be computed in polynomial time.
\end{lemma}

\begin{proof}
  Applying the Cauchy Schwarz inequality, if \(f, g \in \End(A)\), we may compute a good representation of \(f \circ g\) in polynomial time.

  By \cite[Lemma 3.2]{anni2025computation}, we may compute the products \(\langle f, g\rangle_A\), where \(f\) and \(g\) are either elements of \(B\) or products of elements of \(B\).
  Set \(B = (b_i)_{1 \leq i \leq 16}\) and fix \(f, g \in B\), we may then compute the Gram matrix \(G\) of the family \(B \cup f \circ g\). Since \(B\) forms a basis of \(\End(A)\), we may then recover the coordinates of \(f \circ g\) with respect to \(B\). Indeed, let \(v = (v_i)_{1 \leq i \leq 17} \in \Ker(G)\). We have
  \[0 = v^t G v = \langle v_{17} f \circ g + \sum_{i = 1}^{16} v_i b_i,f \circ g + \sum_{i = 1}^{16} v_i b_i\rangle_A.\]
  It follows that \(f \circ g = \sum_{i = 1}^{16} \frac{v_i}{v_{17}} b_i\). We may thus recover the coordinates of the \(b_i \circ b_j\) with respect to \(B\), i.e. the structure constants of \(\End(A)\) with respect to \(B\).
\end{proof}

\begin{lemma}\label{lemma:endo_curves_from_endo_product}
    There is a polynomial-time algorithm which, given supersingular elliptic curves \(E_1\) and \(E_2\), as well as an efficient representation of \(\End(A)\) where \(A = E_1 \times E_2\), outputs efficient representations of \(\End(E_1)\) and \(\End(E_2)\).
\end{lemma}

\begin{proof}
We only detail the computation of \(\End(E_1)\), as that of \(\End(E_2)\) is similar. Let \(f_1, \hdots, f_{16}\) be the basis of \(\End(A)\) such that efficient representations of the \(f_i\) are known. Then, by restricting to \(E_1 \times \{\infty_2\}\) and composing with the projection map \(A \to E_1\), we may deduce efficient representations of \(16\) endomorphisms of \(E_1\) which generate its endomorphism ring as a \(\ZZ\)-module.
Applying \cite[Lemma 8.2]{wesolowski2022supersingular}, we may compute the \(\langle f_i, f_j\rangle_{E_1}\), and therefore a Gram matrix of any subfamily of \((f_1,\hdots,f_{16})\). We may therefore test that a subfamily is free and extract a basis of \(\End(E_1)\).
\end{proof}

%% file: sections/problems.tex
\section{The computational problems}\label{sec:problems}

\begin{problem}[Explicit Ibukiyama-Katsura-Oort correspondence]\label{prob:IKO}
    Given a principally polarised superspecial abelian surface \(A\), compute its Ibukiyama-Katsura-Oort representative, \(\mu(A)\).
\end{problem}

\begin{problem}[Explicit isomorphism problem]\label{prob:Iso}
    Given two principally polarised superspecial abelian surface \(A_1\) and \(A_2\), compute an unpolarised isomorphism \(\varphi\colon A_1 \to A_2\) and its inverse \(\varphi^{-1}\).
\end{problem}

\begin{problem}[Good Endomorphism ring computation]\label{prob:Endo_Computation}
    Given a principally polarised superspecial abelian surface \(A\), compute a good representation of its endomorphism ring.
\end{problem}

\begin{problem}[Effective endomorphism ring computation]\label{prob:Endo_Weak_Computation}
    Given a principally polarised superspecial abelian surface \(A\), compute a weak representation of its endomorphism ring.
\end{problem}

\begin{theorem}\label{thm:main}
  Let \((A, \lambda_A)\) be a PPSAS. Under GRH and the assumptions for applying the \(\mathsf{KLPT}^2\) algorithm, \Cref{prob:IKO,prob:Iso,prob:Endo_Computation,prob:Endo_Weak_Computation} are equivalent if \(A\) is assumed to be a  product of supersingular elliptic curves (and if \(A_1\) and \(A_2\) are assumed to be products for \Cref{prob:Iso}). In the general case, \Cref{prob:IKO} and \Cref{prob:Endo_Computation} are equivalent, and furthermore, \Cref{prob:Endo_Weak_Computation} reduces to \Cref{prob:Iso} and \Cref{prob:Iso} reduces to \Cref{prob:IKO}.
\end{theorem}

\begin{proof}
    The various reductions are the object of the remainder of this article. We provide diagrams indicating the location of each reduction. In these diagrams, arrows point from a problem to the problem it reduces to.\\
    In the case where the varieties are assumed to be products of supersingular elliptic curves, the reductions are as follows:
    \begin{center}
    \begin{tikzpicture}
        \node[draw, align=center] (IKO) at (0, 0) {\Cref{prob:IKO} \\ Explicit IKO};
        \node[draw, align=center] (ISO) at (8, 0) {\Cref{prob:Iso} \\ Explicit isomorphism};
        \node[draw, align=center] (ENDST) at (0, -3) {\Cref{prob:Endo_Computation} \\  Good endomorphism ring};
        \node[draw, align=center] (END) at (8, -3) {\Cref{prob:Endo_Weak_Computation} \\ Effective endomorphism ring};

        \draw[->] (ISO) to node[midway, above=0.1] {\cite[Theorem 3.2]{gaudry2025computing}} (IKO);
        \draw[->] (END) to node[midway, right=0.1] {\Cref{prop:reduction_endo_to_iso}} (ISO);
        \draw[->] (IKO) to node[midway, left=0.1] {\Cref{thm:IKO_to_Endo}}  (ENDST);
        \draw[->] (ENDST) to node[midway, below=0.1] {\Cref{prop:product_known_endo_to_weakly_known}} (END);
    \end{tikzpicture}
\end{center}

  For the general case, we provide the following reductions:
    \begin{center}
    \begin{tikzpicture}
        \node[draw, align=center] (IKO) at (0, 0) {\Cref{prob:IKO} \\ Explicit IKO};
        \node[draw, align=center] (ISO) at (7, 0) {\Cref{prob:Iso} \\ Explicit isomorphism};
        \node[draw, align=center] (ENDST) at (0, -3) {\Cref{prob:Endo_Computation} \\ Good endomorphism ring};
        \node[draw, align=center] (END) at (7, -3) {\Cref{prob:Endo_Weak_Computation} \\ Effective endomorphism ring};

        \draw[->] (ISO) to node[midway, above=0.1] {\Cref{thm:ISO_to_IKO}} (IKO);
        \draw[->] (END) to node[midway, right=0.1] {\Cref{prop:reduction_endo_to_iso}} (ISO);
        \draw[->] (IKO) to[bend right = 30] node[midway, left=0.1] {\Cref{thm:IKO_to_Endo}}  (ENDST);
        \draw[->] (ENDST) to[bend right=30] node[midway, right=0.1] {\Cref{prop:Good_Endo_To_IKO}} (IKO);
    \end{tikzpicture}
\end{center}
\end{proof}

\begin{proposition}\label{prop:reduction_endo_to_iso}
    \Cref{prob:Endo_Weak_Computation} reduces to \Cref{prob:Iso}.
\end{proposition}

\begin{proof}
  Let \(A\) be a given superspecial abelian surface. Let \(\varphi\colon A_0 \to A\) be an isomorphism. Then from an efficient representation of \(\varphi\) and efficient representations of endomorphisms \(f_1,\hdots, f_{16}\) forming a basis of \(\End(A_0)\), we may immediately deduce efficient representations of endomorphisms \(\varphi_*(f_1),\hdots,\varphi_*(f_{16})\), which form a basis of \(\End(A)\).
\end{proof}

%% file: sections/endomorphism_rings.tex
\section{Reductions of Endomorphism Ring Computation}\label{sec:endo}

We treat here two problems regarding the knowledge of the endomorphism ring of a superspecial abelian surface. If the surface is a product, we prove that good endomorphism ring computation reduces to efficient endomorphism ring computation. If the surface is a Jacobian, we instead prove a reduction to the Explicit Ibukiyama-Katsura-Oort correspondence.

\subsection{Endomorphism Ring Computation for a product of supersingular elliptic curves}
Let \(E_1, E_2\) be supersingular elliptic curves.

For the abelian surface \(A = E_1 \times E_2\) with its product polarisation. We will reduce the problem to the case of elliptic curves and then apply \Cref{lemma:curve_knowing_endo}.

\begin{proposition}\label{prop:product_known_endo_to_weakly_known}
  Under GRH, if the endomorphism ring of \(A\) is weakly known, then a good representation of it may be computed in polynomial time.
\end{proposition}
\begin{proof}
  Using \Cref{lemma:endo_curves_from_endo_product}, efficient representations of the endomorphism rings of \(E_1\) and \(E_2\) may be computed in polynomial time. Then, using \Cref{lemma:curve_knowing_endo}, good representations of these endomorphism rings may be computed.
  It follows from the discussion in \Cref{ex:isogeny_matrix} that if the endomorphism rings of \(E_1\) and \(E_2\) are known, then so is the endomorphism ring of \(A\).
\end{proof}

\subsection{Endomorphism Ring Computation from the Ibukiyama-Katsura-Oort matrix}
\begin{proposition}\label{prop:Good_Endo_To_IKO}
    \Cref{prob:Endo_Computation} reduces to \Cref{prob:IKO}.
  \end{proposition}
  \begin{proof}
    Let \(A\) be a PPSAS, whose Ibukiyama-Katsura-Oort matrix \(\mu(A)\) is knwon. By \Cref{thm:ISO_to_IKO}, we may compute an isomorphism \(\varphi\colon A_0 \to A\).
    Using the proof of \Cref{prop:reduction_endo_to_iso}, we may compute a basis \(B\) of \(\End(A)\) and efficient representations of its element, given as \(\varphi_*(b_i)\), where the \(b_i\) form a basis of \(\End(A_0)\). The Rosati involution corresponding to \(\mu(A)\) may be computed using \Cref{lemma:ros_iso_IKO}, and then carried to \(\End(A)\). Likewise, upperbounds for the \(\langle \theta(b_i), \theta(b_i)\rangle_{A}\) may be computed. Since \(\varphi\) induces an isomorphism between \(H^0(A, \Omega_A)\) and \(H^0(A_0, \Omega_{A_0})\), we may also use the matrices \(b_i^*\) to complete the good representations of the elements of \(B\).
  \end{proof}

  \begin{remark}
    This only guarantees that the Rosati involution computed on \(\End(A)\) is induced by a principal polarisation which is merely equivalent to the natural polarisation on \(A\) as introduced in \Cref{ex:natural_principal_representations}. This is not a problem, since we only represent PPSASs up to polarised isomorphism.
  \end{remark}

%% file: sections/isomorphism.tex
\section{Computing unpolarised isomorphisms}\label{sec:iso}

In this section, we aim to prove that \Cref{prob:Iso} reduces to \Cref{prob:IKO}. That is, we shall prove that if \(A_1, A_2\) are PPSASs with known Ibukiyama-Katsura-Oort matrices, an isomorphism \(\varphi\colon A_1 \to A_2\) and its inverse can be computed in polynomial time.

The case where \(A_1\) and \(A_2\) are products of elliptic curves is treated in \cite{gaudry2025computing}. The fact that an isomorphism may easily be inverted in this case is not stated, but follows directly from \Cref{lemma:invert_isom_product}.

In order to complete this result, it is enough to provide a polynomial algorithm which, given a superspecial Jacobian surface \(A\) and the matrix \(\mu(A)\), outputs two superspecial elliptic curves \(E_1\) and \(E_2\), their endomorphism rings, and an isomorphism \(\varphi\colon E_1 \times E_2 \to A\). Our strategy is inspired from the technique used in proving \cite[Theorem 2]{oort1975which} (a superspecial Jacobian surface is isomorphic to a product of supersingular elliptic curves), although some work is required to make the proof computational.

Knowing the matrix \(\mu(A)\), we may compute a chain of \((2,2)\)-isogenies from \(E_0 \times E_0\) to \(A\). We shall work the chain one link at a time, and so we state the following lemma.

\begin{lemma}\label{lemma:iso_comput_2_2}
    Let \(E_1, E_2\) be elliptic curves with known endomorphism rings. Let \(f\colon E_1 \times E_2 \to A\) be a \((2,2)\)-isogeny with known kernel \(K\). Then, assuming GRH, a pair of elliptic curves \(E_3, E_4\), a \((2,2)\)-isogeny \(E_1 \times E_2 \to E_3 \times E_4\) and an unpolarised isomorphism \(E_3 \times E_4 \to A\) can be computed in polynomial time.
\end{lemma}

\begin{proof}
    Let \(P = (P_1, P_2)\) and \(Q = (Q_1, Q_2)\) be generators of \(K\). Since the endomorphism rings of \(E_1\) and \(E_2\) are known, we may compute an isogeny \(f\colon E_1 \to E_2\) of odd degree \cite{wesolowski2022supersingular}.
    Then, the isogenies \(f\) and \(\widehat{f}\) act bijectively on the two-torsion sets of \(E_1\) and \(E_2\).

    Since the ring of endomorphisms of an elliptic curve surjects onto \(M_2(\FF_2)\) when considering the action of endomorphisms on the two-torsion, it is straightforward to compute endomorphisms \(\alpha_1 \in \End(E_1)\) and \(\alpha_2(\End(E_2))\) such that 
    \begin{align*}
        \alpha_2 \circ f(P_1) &= P_2, \\
        \alpha_2 \circ f(Q_1) &= P_2 + Q_2, \\
        \alpha_1 \circ \widehat{f}(P_2) &= Q_1.
    \end{align*}
    We set \(\varphi = \alpha_2 \circ f\) and \(\psi = \alpha_1 \circ \widehat{f}\).

    We may then consider the unpolarised endomorphism
    \[\theta = \begin{pmatrix}
        \Id_{E_1} + \psi \circ \varphi  & \psi \\
        \varphi                         & \Id_{E_2}
    \end{pmatrix}.\]
    A straightforward computation shows that \(\theta\) is an unpolarised automorphism with inverse
    \[\theta^{-1} = \begin{pmatrix}
        \Id_{E_1} & -\psi \\
        -\varphi  & \Id_{E_2} + \varphi \circ \psi
    \end{pmatrix}\]
     and such that \(\theta(K)\) is the group generated by \((P_1, 0)\) and \((0, P_2)\).

     Setting \(E_3 = E_1 / \langle P_1 \rangle\) and \(E_4 = E_2 / \langle P_2 \rangle\), we may explicitely compute the isogeny \(g\colon E_1 \times E_2 \to E_3 \times E_4\) with kernel \(\theta(K)\). Furthermore, since \(\Ker(g) = \Ker(f \circ \theta^{-1})\), there exist an unpolarised isogeny \(\tau\colon E_3 \times E_4\) such that the following diagram commutes:
     
     We get the following diagram:
     \begin{center}\begin{tikzcd}
        E_1 \times E_2 \ar[d, "\theta"] \ar[r, "f"] & A \\
        E_1 \times E_2 \ar[r, "g"] & E_3 \times E_4 \ar[u, "\tau"]
     \end{tikzcd}\end{center}
    Since \(\deg(\theta) = 1\), and \(\deg(f) = \deg(g) = 4\), it follows that \(\tau\) is an unpolarised isomorphism.
    
    If \(P \in E_3 \times E_4\), \(\tau(P)\) may be computed as \(f \circ \theta^{-1}(Q)\), where \(Q\) is any preimage of \(P\) by \(g\), which may itself be computed using \Cref{lemma:preimage_2_2_isogeny}. The inverse isomorphism \(\tau^{-1}\) may be computed in the same manner.
\end{proof}

\begin{algorithm}
    \KwIn{A superspecial Jacobian surface \(A\)}
    \KwIn{Its Ibukiyama-Katsura-Oort matrix \(\mu(A)\)}
    \KwOut{Supersingular elliptic curves \(E_3, E_4\)}
    \KwOut{An unpolarised isomorphism \(\tau\colon E_3 \times E_4 \to A\)}
    \KwOut{\(\tau^{-1}\)}
    \KwOut{The endomorphism rings of \(E_3\) and \(E_4\)}
    \((A_0 \xrightarrow{f_0} A_1 \xrightarrow{f1} \dots \xrightarrow{f_e} A_e) \gets\) Chain of \((2,2)\) isogenies from \(E_0^2\) to \(A\) \nllabel{line:klpt2}\;
    \(f \gets f_0\)\;
    \(E \gets E_0\)\;
    \(E' \gets E_0\)\;
    \For{\(1 \leq i \leq e\)}{
        \tcc{
            Invariant: \(f\) is a \((2,2)\)-isogeny from \(E \times E'\) to \(A_i\) and the endomorphism rings of \(E\) and \(E'\) are known. \\
            Apply \Cref{lemma:iso_comput_2_2} to \(f\):
        }
        \(E, E' \gets\) supersingular elliptic curves\;
        \(\tau \gets E \times E' \simeq A_i\)\;
        \(R_E, R_{E'} \gets \End(E), \End(E')\)\;
        \(f \gets f_i \circ \tau^{-1}\)\;
    }
    \KwRet{\((E, E', \tau, \tau^{-1}, R_E, R_{E'})\)}\;
    \caption{Unpolarised isomorphism from product to Jacobian}
    \label{algo:iso_prod_jac}
\end{algorithm}

\begin{proposition}
    \Cref{algo:iso_prod_jac} is correct and runs in polynomial time.
\end{proposition}

\begin{proof}
    The computation in Line \ref{line:klpt2} is achieved in polynomial time by applying \(\mathsf{KLPT}^2\) as discussed in \Cref{lemma:klpt2_2_2}, using the fact that \(\mu(E_0^2)\) is the identity matrix. Likewise, the number \(e\) of links in the chain is polynomial in the input. The loop invariant is satisfied by \Cref{lemma:iso_comput_2_2}, and the invariant itself ensures both that the algorithm of \Cref{lemma:iso_comput_2_2} can be applied within the loop, and that the output of the algorithm is correct, considering that \(A_e = A\). Finally, each iteration of the loop runs in polynomial time by the same lemma.
\end{proof}

\begin{theorem} \label{thm:ISO_to_IKO}
    Under the assumptions for applying \(\mathsf{KLPT}^2\), \Cref{prob:Iso} reduces to \Cref{prob:IKO}.
\end{theorem}

\begin{proof}
  Let \(A_1, A_2\) be PPSASs, with known Ibukiyama-Katsura-Oort matrices \(\mu(A_1)\) and \(\mu(A_2)\). If \(A_i = E^{(i)}_1 \times E^{(i)}_2\) is a product of elliptic curves, we compute \(\End(E^{(i)}_1)\) and \(\End(E^{(i)}_2)\) using \Cref{lemma:endo_curves_from_endo_product}. Otherwise, we compute supersingular \(E^{(i)}_1\) and \(E^{(i)}_2\), their endomorphism rings, and an unpolarised isomorphism \(E^{(i)}_1 \times E^{(i)}_2 \to A_i\) using \Cref{algo:iso_prod_jac}.

    We may then apply \cite[Theorem 3.2]{gaudry2025computing} to compute an isomorphism from \(E^{(1)}_1 \times E^{(1)}_2 \to E^{(2)}_1 \times E^{(2)}_2\), whose inverse may be computed using \Cref{lemma:invert_isom_product}.
    
    Putting everything together, this solves \Cref{prob:Iso}.
\end{proof}

%% file: sections/endo_to_iko.tex
\section{Recovering the Ibukiyama-Katsura-Oort matrix from the endomorphism ring}\label{sec:endo_to_iko}

In this section, we show how the Ibukiyama-Katsura-Oort matrix of a principally polarised superspecial abelian surface \(A\) may be computed from the knowledge of its endomorphism ring. Essentially, the strategy is to compute an isomorphism from \(\End(A)\) to \(M_2(O_0)\), transport the Rosati involution from \(\End(A)\) to \(M_2(O_0)\), and use it to recover the Ibukiyama-Katsura-Oort matrix. 
Difficulties arise as an isomorphism between endomorphism rings of superspecial abelian surfaces is not necessarily induced by an unpolarised isomorphism between the surfaces, which may lead to the computation of the wrong Rosati involution and therefore the wrong matrix. Fortunately, a "wrong" isomorphism of endomorphism rings may be identified by checking whether its action modulo \(p\) respects the action of endomorphisms on differentials, and it may then be corrected by post-composition with a fixed automorphism of \(A_0\).

\subsection{Isomorphisms of varieties and isomorphisms of endomorphism rings}

If \(f\colon A \to A'\) is an unpolarised isomorhism of abelian varieties, it induces an isomorphim \(f_*\) from \(\End(A)\) to \(\End(A')\). We show below that an isomorphism \(\End(A) \to \End(A')\) is of the form \(f_*\) if and only if it respects the so-called orientations induced by \(A\) on \(\End(A)\) and \(A'\) on \(\End(A')\). This is a generalisation of results on endomorphisms on supersingular elliptic curves exposed in \cite[Section 42.4]{voight2021quaternion} and ultimately credited to Mestre and Oesterlé.

\begin{remark}
    While the term \emph{orientation} is widely used in the field of isogeny-based cryptography to mean an embedding of a quadratic number ring into a quaternion order, it refers here to a widely different notion and there should be no confusion between the two. We keep the term "orientation" so that our exposition mirrors that of \cite[Section 42.4]{voight2021quaternion} better.
\end{remark}

\begin{remark}
    The focus of our work is superspecial abelian surfaces. Most of the theory of orientations that we develop below is however valid in arbitrary dimension greater than \(1\). For the sake of simplicity, we state the results for surfaces only.
\end{remark}

We clarify this strategy with the following lemma:
\begin{lemma}\label{lemma:iso_inner_if_induced}
  Let \(A, A'\) be abelian varieties such that there exists an unpolarised isomorphism \(\varphi_0\colon A' \to A\). Then, for any unpolarised isomorphism, \(\varphi\colon A \to A'\), the induced automorphism \((\varphi_0)_* \circ \varphi_*\) of \(\End(A)\) is inner. 
\end{lemma}

\begin{proof}
  This is clear as \((\varphi_0)_* \circ \varphi_* = (\varphi_0 \circ \varphi)_*\) is the inner automorphism of \(\End(A)\) induced by \(\varphi_0 \circ \varphi\).
\end{proof}

Below, we prove that the converse is true. We also show that the group of outer automorphisms of \(\End(A)\) has order two, and we provide an explicit criterion for an isomorphism of endomorphism rings to be induced by an isomorphism of abelian varieties.

\subsubsection{Orientations}
\begin{definition}\label{def:orientation}
  Let \(R\) be a ring isomorphic to a maximal order in \(M_2(B_{p,\infty})\). An \emph{orientation} of \(R\) is a surjective homomorphism \(\zeta\colon R \to M_2(\overline{\FF_{q}})\) such that \(\Ker(\zeta) = P\), where \(P\) is the intersection of \(R\) and the unique two-sided ideal of \(R_p\) (see \cite[Theorem 18.3]{reiner1975maximal}).

    Two orientations \(\zeta\) and \(\zeta'\) are \emph{equivalent}, denoted by \(\zeta \sim \zeta'\), if \(\det \circ \zeta = \det \circ \zeta'\).

    An \emph{oriented maximal order} is the data of a ring \(R\) as above, together with an equivalence class of orientations. An \emph{isomorphism of oriented maximal orders} \(\varphi\colon (R, [\zeta]) \to (R', [\zeta'])\) is then a ring isomorphism from \(R\) to \(R'\) such that \(\zeta \sim \zeta' \circ \varphi\). We also say that a ring isomorphism is \emph{oriented} if it is an isomorphism of oriented maximal orders, when the orientations are clear from context.
\end{definition}

The facts we state below about the ideal \(P\) are stated in \cite[Chapter 5]{reiner1975maximal}.

Let \(R\) be as in \Cref{def:orientation}. There are exactly two classes of orientations on \(R\), and they are related in a straightforward manner by any outer automorphism of \(R\), which we describe in the two lemmas below.

We write \(\Aut(R)\) for the group of automorphisms of \(R\), \(\Inn(R)\) for the group of inner automorphisms of \(R\) and \(\Out(R)\) for the group \(\Aut(R)/\Inn(R)\) of \emph{outer automorphisms} of \(R\).

\begin{lemma}\label{lemma:out_R}
    The group \(\Out(R)\) is isomorphic to \(\ZZ/2\ZZ\). Its non-trivial class contains, in particular, any automorphism defined as conjugation by a generator of \(P\) as a left ideal.
\end{lemma}

\begin{proof}
  Since \(M_2(B_{p, \infty})\) has discriminant \(p^2\), \cite[Corollary 37.32]{reiner1975maximal} implies that the Picard group of \(R\) has order \(2\), and is generated by \(P\) as a two-sided ideal. Since the class group of \(R\) is trivial \cite[Theorem 2.1]{ibukiyama1986supersingular}, it then follows from \cite[Corollary 37.34]{reiner1975maximal} that \(\Out(R) \simeq \Pic(R)\) and therefore has order \(2\). To avoid any confusion, we direct the attention of the reader to the fact that the ring that we call \(R\) is denoted by \(\Lambda\) in \cite{reiner1975maximal}, while the ring denoted there by \(R\) is simply \(\ZZ\) in our case.

  Now, since all maximal orders in \(M_2(B_{p, \infty})\) are isomorphic (this follows from \cite[Theorem 21.6 and Theorem 27.8]{reiner1975maximal} together with the fact that the class group of \(R\) is trivial), we may assume that \(R = M_2(O)\), where \(O\) is a maximal order of \(M_2(B_{p, \infty})\) containing some \(j\) such that \(\nrd(j) = p\). Then, by the discussion aboce \cite[Corollary 42.4.12]{voight2021quaternion}, conjugation by \(j\) induces an outer automorphism of \(O\), and therefore conjugation by \(j I_2\) induces an outer automorphism of \(R\). It is clear in this case that \(j I_2\) is a generator of \(P\) as a left ideal.
\end{proof}

\begin{lemma}\label{lemma:classes_orientations}
    Let \(R\) be as in \Cref{def:orientation}. There are exactly \(2\) equivalence classes of orientations on \(R\). If \(f \in \Aut(R)\) and \(\zeta\) is an orientation of \(R\), then \(\zeta \sim \zeta \circ f\) if and only if \(f \in \Inn(R)\).
\end{lemma}

\begin{proof}
    By definition, an orientation \(\zeta\) factors through the map \(\pi\colon R \to R/P \simeq M_2(\FF_{q})\). Since an orientation is surjective, it is fully described by the choice of an automorphism \(\varphi\) of \(M_2(\FF_q)\), through the relation \(\zeta = \varphi \circ \pi\). Such an automorphism \(\varphi\) factors uniquely as \(\iota \circ \sigma^\varepsilon\), where \(\iota \in \Inn(M_2(\FF_{p^2}))\), \(\sigma\) is induced by the non-trivial automorphism of \(\FF_{p^2}\) and \(\varepsilon \in \{0, 1\}\). It is clear that two orientations are equivalent if and only if the \(\varepsilon\) variables as they appear in their representations as above are equal. 

    If \(u \in R^\times\), and \(\iota_u\) the inner automorphism of \(R\) induced by \(u\), we have \(\zeta \circ \iota_u\colon a \mapsto \zeta(u) \zeta(a) \zeta(u)^{-1}\), and therefore \(\zeta \sim \zeta \circ \iota_u\). Conversely, if \(f \in \Aut(R) \setminus \Inn(R)\), we identify \(R\) with \(M_2(O)\) as in the proof of \Cref{lemma:out_R}. Then, we may assume that \(f\) is conjugation by \(j I_2\), which induces the Galois conjugation on \(R / P \simeq M_2(\FF_{p^2})\). It is then clear that \(\zeta \not\sim \zeta \circ f\).
\end{proof}

\subsubsection{Orientations and superspecial abelian surfaces}

Let \(A\) be a superspecial abelian surface, so that \(\End(A)\) is isomorphic to a maximal order in \(M_2(B_{p, \infty})\). We define an orientation on \(\End(A)\) induced by \(A\) using the action of endomorphisms of \(A\) on its space of global differentials \(H^0(A, \Omega_A)\).

Recall that the space \(D(A) \coloneqq H^0(A, \Omega_{A/\FF_q})\) is a \(2\) dimensional \(\FF_q\)-vector space. An endomorphism of \(A\) acts on \(D(A)\) as a pullback map \(f^*\). It follows that a choice of basis of \(H^0(A, \Omega_A)\) yields a map
\[\begin{array}{rccc}
    \zeta_A\colon &\End(A) &\to &M_2(\FF_q) \\
    & f &\mapsto &f^*.
\end{array}\]

The map \(\zeta_A\) is well defined up to post-composition by an inner automorphism of \(M_2(\FF_q)\), which is enough to give \(\End(A)\) a natural structure of oriented maximal order, after we prove the following lemma.

\begin{lemma}
    The map \(\zeta_A\) defined above induces an orientation of the order \(R_A\) identified with \(\End(A)\).
\end{lemma}

\begin{proof}
    We must check that the ring homomorphism \(\zeta_A\) is surjective, and that it kills endomorphisms with degree divisible by \(p\).

    Since polarisations are uninvolved in this result, we may assume that \(A = E^2\) is a product of supersingular elliptic curves. Then we may identify \(\End(A)\) with \(M_2(\End(E))\), and \(D(A) = D(E)^2\). Now, let \(\zeta_E\) be the map \(\End(E) \to \FF_{q}\) induced by the action of \(\End(E)\) on \(D(E)\). It is clear that the matrix \((f_{ij})_{1 \leq i,j \leq 2} \in M_2(\End(E))\) acts on \(D(A)\) as the matrix \((f^*_{ij})_{1 \leq i, j \leq 2}\). In particular, the map \(\zeta_A\) is surjective.

    Furthermore, if an endomorphism of \(\End(A)\) lying in \(P\), is written as a matrix of isogenies, each of its component is inseparable. Indeed, they all have degree divisible by \(p\) and are therefore inseparable because the curve \(E\) is supersingular. Then, it is known that such isogenies act as zero on the global differentials of \(E\) (see \cite[Corollary 5.6.(c)]{silverman2009arithmetic}), and therefore \(f\) as a whole acts as \(0\) on all global differentials of \(A\).
\end{proof}

The point of introducing orientations is that, as in the genus \(1\) case, isomorphisms of varieties induce isomorphisms of oriented endomorphism rings:

\begin{theorem}\label{thm:iso_endo_var}
    Let \(A, A'\) be superspecial abelian surfaces, and let \(\tau\colon \End(A) \to \End(A')\) be an isomorphism of rings.

    Then there exists some unpolarised isomorphism \(\varphi\colon A \to A'\) such that \(\tau = \varphi_*\) if and only if \(\tau\) is a morphism of oriented maximal orders from \(\End(A), \zeta_A\) to \(\End(A'), \zeta_{A'}\).
\end{theorem}

\begin{proof}
    Let \(\varphi\colon A \to A'\) be an isomorphism, and consider the isomorphism \(\varphi^* \colon D(A') \to D(A)\). It is clear from the definitions that if \(f \in \End(A)\), then 
    \[\zeta_{A'} \circ \varphi_* (f) = (\varphi^{-1})^* \circ \zeta_A(f) \circ \varphi^*.\]
    It follows that the isomorphism \(\varphi_*\) is oriented as required.

    Let \(\tau\colon \End(A) \to \End(A')\) be a non-oriented isomorphism. Now, fix an isomorphism \(\varphi_0\colon A' \to A\) and set \(\tau_0 = (\varphi_0)_*\) the induced morphism of oriented maximal orders. Since \(\tau_0\) is oriented, it follows that \(\tau_0 \circ \tau\) is not, as clearly the composition of two oriented isomorphisms is itself oriented. By \Cref{lemma:classes_orientations}, it follows that \(\tau_0 \circ \tau\) is not an inner automorhism of \(\End(A)\). By \Cref{lemma:iso_inner_if_induced}, it follows that \(\tau\) is not induced by any isomorphism of abelian varieties.

\end{proof}

\subsection{The reduction}

We may now state the reduction from \Cref{prob:IKO} to \Cref{prob:Endo_Computation}. We let \(\zeta_0\) be the orientation of \(M_2(O_0)\) induced by \(A_0\).  We also let \(\pi \in M_2(O_0)\) be a generator of the right ideal \(P_0 = \{a \in M_2(O_0)\ \colon\ p \mid \nrd(a)\}\).

\begin{algorithm}
    \KwIn{A superspecial abelian surface \(A\) with known endomorphism ring}
    \KwOut{The Ibukiyama-Katsura-Oort matrix \(\mu(A)\)}
    \(\tau \gets\) isomorphism from \(\End^0(A)\) to \(M_2(B_{p, \infty})\) \nllabel{line:compute_explicit_isom}\;
    \(R \gets \tau(\End(A))\) \nllabel{line:apply_tau}\;
    \(\gamma \gets \) generator of the right \(M_2(O_0)\)-ideal \(R M_2(O_0)\) \nllabel{line:solve_PIP}\;
    \If{\(\exists b \in \End(A)\ \mid\ \det(\zeta_A(b)) \neq \det(\zeta_{A_0}(\gamma^{-1} \tau(b) \gamma))\)}{ \nllabel{line:if}
        \(\gamma \gets \gamma \pi\) \nllabel{line:multiply}\;
    }\nllabel{line:end_if}
    \(\sigma \gets \left(\begin{array}{ccc}
                    M_2(B_{p, \infty}) &\to & M_2(B_{p, \infty}) \\
                    \alpha &\mapsto &\gamma^{-1}(\tau \circ r_A \circ \tau^{-1})(\gamma \alpha \gamma^{-1}) \gamma
                 \end{array}\right)\) \nllabel{line:compute_sigma} \;
    Compute \(g \in Mat(A_0)\) corresponding to the involution \(\sigma\)\; \nllabel{line:compute_g_from_rosati}
    \KwRet{\(g\)}\;
    \caption{Computing the Ibukiyama-Katsura-Oort matrix from the endomorphism ring}
    \label{algo:IKO_to_Endo}
\end{algorithm}

\begin{theorem}\label{thm:IKO_to_Endo}
    Assuming GRH, \Cref{algo:IKO_to_Endo} is correct and runs in polynomial time.
\end{theorem}

\begin{proof}
We first prove that the theorem runs in polynomial time. Since there are no loops or recursive calls, we only need to check that each line runs in polynomial time. Lines \ref{line:apply_tau}, \ref{line:multiply} and \ref{line:compute_sigma} are straightforward computations involving linear algebra over \(\QQ\). As the inputs for these computations are either part of the input of the algorithm or outputs of other lines, these lines may be executed in polynomial time if the rest of the algorithm does as well. The problem in Line \ref{line:compute_explicit_isom} may be solved using \Cref{lemma:structure_constants,lemma:explicit_isom_CSA}. The problem in Line \ref{line:solve_PIP} is an instance of a principal ideal problem, which is solved by \cite[Theorem 3.13]{page2026conversion}. The only difficulty in Line \ref{line:if} would be to compute \(\zeta_A\), but it is known since good representations of the basis elements of \(\End(A)\) are known. Finally, the computation in Line \ref{line:compute_g_from_rosati} may be done using \Cref{lemma:compute_mu_from_rosati}.

    Now, the ideal \(I = RM_2(O_0)\) is both a right \(M_2(O_0)\)-ideal and a left \(R\)-ideal. If \(\gamma\) is a generator of \(I\) as a right \(M_2(O_0)\)-ideal, then for any \(a \in R\) and \(b \in M_2(O_0)\), we have \(a \gamma b = \gamma c\) for some \(c \in M_2(O_0)\). Since \(R\) is a maximal order, it follows that \(R = \gamma M_2(O_0) \gamma^{-1}\), and that conjugation by \(\gamma\) induces an isomorphism between the two rings.

    Therefore, the map \(b \mapsto \gamma^{-1} \tau(b) \gamma\) is an isomorphism from \(\End(A)\) to \(M_2(O_0) = \End(A_0)\). It then follows from \Cref{lemma:out_R,lemma:classes_orientations} that after Line \ref{line:end_if}, the map \(b \mapsto \gamma^{-1} \tau(b) \gamma\) is an oriented isomorphism. By \Cref{thm:iso_endo_var}, there is an unpolarised isomorphism \(\varphi\colon A \to A_0\) such that for all \(b \in \End(A)\), \(\varphi_*(b) = \gamma^{-1} \tau(b) \gamma\). Now, let \(\lambda = \widehat{\varphi^{-1}} \circ \lambda_A \circ \varphi^{-1}\) be the polarisation on \(A_0\) such that \(\varphi\) is a polarised isomorphism from \((A, \lambda_A)\) to \((A_0, \lambda)\). It is straightforward to check that the map \(\sigma\) computed in Line \ref{line:compute_sigma} is in fact the Rosati involution \(r_\lambda\), and that the algorithm therefore outputs \(\mu(\lambda) = \mu(A)\).
\end{proof}

%% file: biblio.bib
@misc{gaudry2025computing,
  title = {Computing Isomorphisms between Products of Supersingular Elliptic
           Curves},
  author = {Pierrick Gaudry and Julien Soumier and Pierre-Jean Spaenlehauer},
  year = {2025},
  eprint = {2503.21535},
  archivePrefix = {arXiv},
  primaryClass = {math.NT},
  url = {https://arxiv.org/abs/2503.21535},
}

@inbook{milne1984jacobian,
  author = {Milne, J. S.},
  title = {Arithmetic Geometry},
  chapter = {Jacobian Varieties},
  pages = {167--212},
  year = {1984},
  publisher = {Springer-Verlag},
  address = {New York},
}

@inbook{milne1984abelian,
  author = {Milne, J. S.},
  title = {Arithmetic Geometry},
  chapter = {Abelian Varieties},
  pages = {103--150},
  year = {1984},
  publisher = {Springer-Verlag},
  address = {New York},
}

@article{igusa1960arithmetic,
  ISSN = {0003486X, 19398980},
  URL = {http://www.jstor.org/stable/1970233},
  author = {Jun-Ichi Igusa},
  journal = {Annals of Mathematics},
  number = {3},
  pages = {612--649},
  publisher = {[Annals of Mathematics, Trustees of Princeton University on
               Behalf of the Annals of Mathematics, Mathematics Department,
               Princeton University]},
  title = {Arithmetic Variety of Moduli for Genus Two},
  urldate = {2026-01-07},
  volume = {72},
  year = {1960},
}

@article{castryck2020hash,
  title = {Hash functions from superspecial genus-2 curves using Richelot
           isogenies},
  author = {Wouter Castryck and Thomas Decru and Benjamin Smith},
  pages = {268--292},
  volume = {14},
  number = {1},
  journal = {Journal of Mathematical Cryptology},
  doi = {doi:10.1515/jmc-2019-0021},
  year = {2020},
}

@inproceedings{wesolowski2022supersingular,
  author = {Wesolowski, Benjamin},
  booktitle = {2021 IEEE 62nd Annual Symposium on Foundations of Computer
               Science (FOCS)},
  title = {The supersingular isogeny path and endomorphism ring problems are
           equivalent},
  year = {2022},
  volume = {},
  number = {},
  pages = {1100-1111},
  doi = {10.1109/FOCS52979.2021.00109},
}

@article{oort1975which,
  author = {Oort, Frans},
  title = {Which abelian surfaces are products of elliptic curves?},
  fjournal = {Mathematische Annalen},
  journal = {Math. Ann.},
  issn = {0025-5831},
  volume = {214},
  pages = {35--47},
  year = {1975},
  language = {English},
  doi = {10.1007/BF01428253},
  url = {https://eudml.org/doc/162679},
  zbMATH = {3443763},
  Zbl = {0283.14007},
}

@inproceedings{castryck2025klpt,
  title = {KLPT 2: algebraic pathfinding in dimension two and applications},
  author = {Castryck, Wouter and Decru, Thomas and Kutas, P{\'e}ter and Laval,
            Abel and Petit, Christophe and Ti, Yan Bo},
  booktitle = {Annual International Cryptology Conference},
  pages = {167--200},
  year = {2025},
  organization = {Springer},
}

@article{csahok2022explicit,
  author = {Csah{\'o}k, T{\'{\i}}mea and Kutas, P{\'e}ter and Montessinos, Micka
            {\"e}l and Z{\'a}br{\'a}di, Gergely},
  title = {Explicit isomorphisms of quaternion algebras over quadratic global
           fields},
  fjournal = {Research in Number Theory},
  journal = {Res. Number Theory},
  issn = {2522-0160},
  volume = {8},
  number = {4},
  pages = {24},
  note = {Id/No 77},
  year = {2022},
  language = {English},
  doi = {10.1007/s40993-022-00380-3},
  zbMATH = {7600669},
  Zbl = {1534.11136},
}

@misc{edixhovenabelian,
  author = {Edixhoven, Bas and van der Geer, Gerard and Moonen, Ben},
  title = {Abelian Varieties},
  howpublished = {\url{http://van-der-geer.nl/~gerard/AV.pdf}},
  note = {Last consulted on 08/01/2026},
}

@article{ibukiyama1986supersingular,
  author = {Ibukiyama, Tomoyoshi and Katsura, Toshiyuki and Oort, Frans},
  title = {Supersingular curves of genus two and class numbers},
  fjournal = {Compositio Mathematica},
  journal = {Compos. Math.},
  issn = {0010-437X},
  volume = {57},
  pages = {127--152},
  year = {1986},
  language = {English},
  url = {https://eudml.org/doc/89752},
  zbMATH = {3946334},
  Zbl = {0589.14028},
}

@book{voight2021quaternion,
  author = {Voight, John},
  title = {Quaternion algebras},
  fseries = {Graduate Texts in Mathematics},
  series = {Grad. Texts Math.},
  issn = {0072-5285},
  volume = {288},
  isbn = {978-3-030-56692-0; 978-3-030-57467-3; 978-3-030-56694-4},
  year = {2021},
  publisher = {Cham: Springer},
  language = {English},
  doi = {10.1007/978-3-030-56694-4},
  zbMATH = {7261776},
  Zbl = {1481.11003},
}

@book{reiner1975maximal,
  author = {Reiner, Irving},
  title = {Maximal orders},
  fseries = {London Mathematical Society Monographs},
  series = {Lond. Math. Soc. Monogr.},
  volume = {5},
  year = {1975},
  publisher = {Academic Press, London},
  language = {English},
  zbMATH = {3477356},
  Zbl = {0305.16001},
}

@misc{page2026conversion,
  author = {Page, Aurel and Robert, Damien and Soumier, Julien},
  title = {On the conversion of module representations for higher dimensional
           supersingular isogenies},
  note = {(in preparation)},
}

@misc{anni2025computation,
  author = {Anni, Samuele and Bisson, Gaetan and Iezzi, Annamaria and Garc{\'{\i
            }}a, Elisa Lorenzo and Wesolowski, Benjamin},
  title = {On the computation of endomorphism rings of abelian surfaces over
           finite fields},
  year = {2025},
  howpublished = {Preprint, {arXiv}:2503.08925 [math.{NT}] (2025)},
  url = {https://arxiv.org/abs/2503.08925},
  arXiv = {arXiv:2503.08925},
}

@article{jordan2018abelian,
  author = {Jordan, Bruce W. and Keeton, Allan G. and Poonen, Bjorn and Rains,
            Eric M. and Shepherd-Barron, Nicholas and Tate, John T.},
  title = {Abelian varieties isogenous to a power of an elliptic curve},
  journal = {Compos. Math.},
  issn = {0010-437X},
  volume = {154},
  number = {5},
  pages = {934--959},
  year = {2018},
  language = {English},
  doi = {10.1112/S0010437X17007990},
}

@inproceedings{eisentrager2018supersingular,
  title = {Supersingular isogeny graphs and endomorphism rings: reductions and
           solutions},
  author = {Eisentr{\"a}ger, Kirsten and Hallgren, Sean and Lauter, Kristin and
            Morrison, Travis and Petit, Christophe},
  booktitle = {Annual International Conference on the Theory and Applications of
               Cryptographic Techniques},
  pages = {329--368},
  year = {2018},
  organization = {Springer},
}

@inproceedings{castryck2023efficient,
  title = {An efficient key recovery attack on SIDH},
  author = {Castryck, Wouter and Decru, Thomas},
  booktitle = {Annual international conference on the theory and applications of
               cryptographic techniques},
  pages = {423--447},
  year = {2023},
  organization = {Springer},
}

@article{kohel2014quaternion,
  title = {On the quaternion-isogeny path problem},
  author = {Kohel, David and Lauter, Kristin and Petit, Christophe and Tignol,
            Jean-Pierre},
  journal = {LMS Journal of Computation and Mathematics},
  volume = {17},
  number = {A},
  pages = {418--432},
  year = {2014},
  publisher = {London Mathematical Society},
}

@misc{robert2024module,
  author = {Damien Robert},
  title = {The module action for isogeny based cryptography},
  howpublished = {Cryptology {ePrint} Archive, Paper 2024/1556},
  year = {2024},
  url = {https://eprint.iacr.org/2024/1556},
}

@book{silverman2009arithmetic,
  author = {Silverman, Joseph H.},
  title = {The arithmetic of elliptic curves},
  edition = {2nd ed.},
  fseries = {Graduate Texts in Mathematics},
  series = {Grad. Texts Math.},
  issn = {0072-5285},
  volume = {106},
  isbn = {978-0-387-09493-9; 978-0-387-09494-6},
  year = {2009},
  publisher = {New York, NY: Springer},
  language = {English},
  doi = {10.1007/978-0-387-09494-6},
  zbMATH = {5549721},
  Zbl = {1194.11005},
}

@inproceedings{maino2023direct,
  title = {A direct key recovery attack on SIDH},
  author = {Maino, Luciano and Martindale, Chloe and Panny, Lorenz and Pope,
            Giacomo and Wesolowski, Benjamin},
  booktitle = {Annual International Conference on the Theory and Applications of
               Cryptographic Techniques},
  pages = {448--471},
  year = {2023},
  organization = {Springer},
}

@inproceedings{robert2023breaking,
  title = {Breaking SIDH in polynomial time},
  author = {Robert, Damien},
  booktitle = {Annual International Conference on the Theory and Applications of
               Cryptographic Techniques},
  pages = {472--503},
  year = {2023},
  organization = {Springer},
}

@inproceedings{flynn2019genus,
  title = {Genus two isogeny cryptography},
  author = {Flynn, E Victor and Ti, Yan Bo},
  booktitle = {International Conference on Post-Quantum Cryptography},
  pages = {286--306},
  year = {2019},
  organization = {Springer},
}

@article{howe2000large,
  author = {Howe, Everett W. and Lepr{\'e}vost, Franck and Poonen, Bjorn},
  title = {Large torsion subgroups of split {Jacobians} of curves of genus two
           or three},
  fjournal = {Forum Mathematicum},
  journal = {Forum Math.},
  issn = {0933-7741},
  volume = {12},
  number = {3},
  pages = {315--364},
  year = {2000},
  language = {English},
  doi = {10.1515/form.2000.008},
  zbMATH = {1425658},
  Zbl = {0983.11037},
}

@misc{herledan2025unconditional,
  author = {Arthur Herlédan Le Merdy and Benjamin Wesolowski},
  title = {Unconditional foundations for supersingular isogeny-based
           cryptography},
  howpublished = {Cryptology {ePrint} Archive, Paper 2025/271},
  year = {2025},
  url = {https://eprint.iacr.org/2025/271},
}
